\documentclass[journal,twoside,web]{ieeecolor}
\usepackage{generic}
\def\BibTeX{{\rm B\kern-.05em{\sc i\kern-.025em b}\kern-.08em
    T\kern-.1667em\lower.7ex\hbox{E}\kern-.125emX}}

\usepackage{textcomp}

\usepackage{sansmath}
\usepackage{graphicx}
\usepackage{amsmath}
\usepackage{xspace}
\usepackage{multirow}
\usepackage{diagbox}
\usepackage{amsfonts}
\usepackage{mathtools}
\usepackage{hyperref}
\usepackage[english]{babel}
\usepackage[version=4]{mhchem}
\usepackage{siunitx}
\usepackage{blindtext}
\usepackage{longtable}
\usepackage{tabularx}
\usepackage[utf8]{inputenc}
\usepackage{subcaption}
\usepackage{algorithm}
\usepackage{algpseudocode}
\usepackage{float}
\usepackage{titlesec}
\usepackage{tikz}
\usepackage{pgfplots}
\usepgfplotslibrary{groupplots}
\pgfplotsset{compat=1.17}
\usepackage{booktabs}
\usepackage{pgfplotstable}
\usepackage{amsthm}
\usepackage{bm}
\usepackage{cancel}
\usepackage{cite}
\usepackage{mathrsfs}
\usepackage{etoolbox}

\newcommand{\ones}{\bm 1}

\DeclareMathOperator*{\Tr}{Tr}

\DeclareMathOperator*{\Cov}{Cov}
\DeclareMathOperator*{\diag}{diag}

\DeclareMathOperator{\blkdiag}{blkdiag}

\newcommand{\mcl}[1]{\mathcal{#1}}
\newcommand{\mbb}[1]{\mathbb{#1}}
\newcommand{\mbf}[1]{\mathbf{#1}}

\DeclareMathOperator*{\minimize}{\mathrm{minimize}}
\DeclareMathOperator*{\subjectto}{\mathrm{subject~to}}

\theoremstyle{definition}
\newtheorem{remark}{Remark}
\newtheorem{proposition}{Proposition}

\newtheorem{definition}{Definition}
\newtheorem{problem}{Problem}
\newtheorem{theorem}{Theorem}
\newtheorem{assumption}{Assumption}

\renewcommand{\baselinestretch}{0.85}

\title{Distributed and Localized Covariance Control of Coupled Systems: A System Level Approach}

\author{Ahmed Khalil, Yoonjae Lee, and Efstathios Bakolas%
\thanks{The first two authors contributed equally as co-first authors.}
\thanks{The authors are with the Department of Aerospace Engineering and Engineering Mechanics, the University of Texas at Austin,  Austin, Texas 78712-1221, USA. Email addresses:
        {\tt\small akhalil@utexas.edu, yol033@utexas.edu, bakolas@austin.utexas.edu}.}
}   

\begin{document}

\maketitle
\thispagestyle{empty}
\pagestyle{empty}

\begin{abstract}
This work is concerned with the finite-horizon optimal covariance steering of networked systems governed by discrete-time stochastic linear dynamics. In contrast to existing work that has only considered systems with dynamically decoupled `agents,' we consider a
dynamically coupled system composed of interconnected
subsystems. In particular, we propose a \emph{distributed} algorithm to compute the \emph{localized} covariance steering policy for each individual subsystem. To this end, we first formulate the localized covariance steering problem and, leveraging the system-level synthesis (SLS) framework, recast this problem as a convex optimization problem over system responses. We then examine the problem's separability and introduce a problem transformation to address instances with nonseparable objective functions. Finally, we develop a distributed algorithm based on the consensus alternating direction method of multipliers (ADMM) to distribute computation across subsystems based on their local information and communication constraints. We demonstrate the effectiveness of our proposed algorithm on a power system with 36 interconnected subsystems.
\end{abstract}


\section{Introduction} \label{sec:introduction}

\IEEEPARstart{T}{he} widespread emergence of large-scale, highly interconnected systems highlights the critical need for scalable controllers that rely only on local information. In complex, safety-critical systems such as power grids, these controllers must ensure the network's reliable operation and optimal performance despite inherent uncertainties and noise. A promising control paradigm addressing these challenges is Covariance Steering (CS), which focuses on guiding the evolution of the system's state distribution while enforcing distributional constraints at the terminal state. By explicitly accounting for stochastic disturbances, CS ensures that the controller achieves its desired mean and covariance targets. CS differs from standard stochastic optimal control methods, such as linear-quadratic-Gaussian (LQG) control \cite{Stengel94}, where uncertainty is only controlled implicitly through shaping the cost function.


Infinite-horizon variants of CS problems were first formulated and studied in the 1980s \cite{Collins85, Hotz1987}. Recently, the finite-horizon counterparts have gained traction, particularly within the controls community \cite{chen2015optimal, bakolas2018finite, Balci22, Balci2021CDC, Balci2021, Okamoto2018, Balci24}. CS has been successfully applied in areas such as path planning \cite{Okamoto19, knaup2023safe, Balci2021}, trajectory optimization \cite{Balci22, Balci2021CDC, Okamoto2018}, and robotic manipulation \cite{lee2022hierarchical}. For discrete-time systems, previous works have presented ways to formulate finite-horizon CS problems as convex optimization problems \cite{liu2024} by utilizing various parameterizations, such as the state history feedback policy \cite{Skaf2010}, the disturbance feedback policy \cite{Balci24}, and the auxiliary variable policy \cite{Okamoto19}.


Distributed CS problems for multi-agent systems have been addressed using similar parametrizations \cite{Saravanos21, Saravanos22, Saravanos23}. However, existing approaches overlook the dynamic coupling and the localized nature of interconnected systems, treating the latter as decoupled `agents.' In tightly coupled systems like power grids, however, distributed controllers require accounting for interactions between subsystems, as disturbances can propagate and amplify across the network, potentially leading
to undesirable effects. Designing such a controller, however, is known to be challenging due to nonconvexity issues that arise when attempting to apply parameterizations \cite{Anderson2019}. This motivates our use of the System Level Synthesis (SLS) framework \cite{wang2018, Anderson2019,Alonso2020DLMPC,Alonso2023}, a method to synthesize a controller over the closed-loop behavior of the system rather than the controller itself. This approach particularly allows for the enforcement of local communication constraints while preserving the problem's convex structure, admitting the use of distributed optimization techniques such as the Alternating Direction Method of Multipliers (ADMM) \cite{Boyd11}.

In this work, we address the finite-horizon \emph{localized} CS problem involving dynamically coupled systems subject to discrete-time stochastic linear dynamics. In sharp contrast to previous work that has utilized the SLS framework to formulate distributed and localized Model Predictive Control (MPC) problems under the case of no driving noise \cite{Alonso2020DLMPC} or deterministic and bounded noise \cite{Alonso2023}, this work focuses on the distributed and localized steering of state distributions in the presence of Gaussian random disturbances. Our contributions are three-fold. First, we show that the localized CS problem can be formulated as a convex optimization problem over system responses via the SLS framework. Second, we introduce a problem transformation to address instances with nonseparable objectives while maintaining the separable structures of certain constraints. Lastly, we propose a distributed algorithm with provable convergence, based on consensus ADMM \cite{Boyd11}, to distribute computation across subsystems while accounting for their local communication and information constraints.


The rest of the paper is structured as follows. Section \ref{sec:probform} introduces the localized covariance steering problem in coupled systems. Section \ref{sec:mainresults} presents a system-level parametrization approach and its distributed solution method. Section \ref{sec:numericalsimulations} contains a numerical simulation of the proposed method. Section \ref{sec:conclusion} concludes this work with directions for future research.

\textit{Notation}: Let $\mbb{N}$ denote the set of non-negative integers, and let $\mbb{S}_{+}^{n}$ and $\mbb{S}_{++}^{n}$ denote the convex cones of $n \times n$ (real symmetric) positive semidefinite and positive definite matrices, respectively. Let $\ones$ denote the vector of ones. The mean and covariance of a random vector $x$ are denoted as $\mbb{E}[x]$ and $\Cov[x]$, respectively. A Gaussian random vector with mean $\mu$ and covariance $\Sigma$ is denoted as $x \sim \mcl{N}(\mu,\Sigma)$. Let $\diag(a_1, \dots, a_N)$ denote a diagonal matrix with scalars $a_1, \dots, a_N$ and $\blkdiag(A_1, \dots, A_N)$ denote a block diagonal matrix with matrices $A_1, \dots, A_N$. The Frobenius norm of a matrix $A$ is denoted as $\|A\|_F \coloneqq \sqrt{\Tr(A^\top A)}$. Let
$A(\mcl{R},\mcl{C})$ denote the submatrix of a matrix $A$ formed by selecting the rows indexed by $\mcl{R}$ and the columns indexed by $\mcl{C}$, with $:$ denoting all indices. Lastly, $\circ$ denotes the elementwise matrix multiplication (Hadamard product).

\section{Problem Formulation} \label{sec:probform}

Consider a discrete-time, linear time-varying (LTV) stochastic system composed of $N$ interconnected subsystems:
\begin{align} \label{eqn:linearsystem}
    \underbrace{\begin{bmatrix}
        x_{t+1}^1
        \\
        \vdots
        \\
        x_{t+1}^N
    \end{bmatrix}}_{x_{t+1}}
    &=
    \underbrace{\begin{bmatrix}
        A_t^{11} & \dots & A_t^{1N}
        \\
        \vdots & \ddots & \vdots
        \\
        A_t^{N1} & \dots & A_t^{NN}
    \end{bmatrix}}_{A_t}
    \underbrace{\begin{bmatrix}
        x_{t}^1
        \\
        \vdots
        \\
        x_{t}^N
    \end{bmatrix}}_{x_t}
    \nonumber
    \\
    &\quad +
    \underbrace{\begin{bmatrix}
        B_t^{1} & &
        \\
        & \ddots &
        \\
        & & B_t^{N}
    \end{bmatrix}}_{B_t}
    \underbrace{\begin{bmatrix}
        u_t^1
        \\
        \vdots
        \\
        u_t^N
    \end{bmatrix}}_{u_t}
    +
    \underbrace{\begin{bmatrix}
        w_t^1
        \\
        \vdots
        \\
        w_t^N
    \end{bmatrix}}_{w_t}
    ,
\end{align}
where $x_t^i \in \mbb{R}^{n_i}$, $u_t^i \in \mbb{R}^{m_i}$, and $w_t^i \in \mbb{R}^{n_i}$ denote the local state, control input, and process noise of subsystem $i$, respectively, with $A_t^{ij} \in \mbb{R}^{n_i \times n_j}$ and $B_t^i \in \mbb{R}^{n_i \times m_i}$. The local disturbance $w_t^i$ is assumed to be a white Gaussian noise process with $\mbb{E}[w_t^i] = 0$ and $\mbb{E}[w_{t'}^i {w_t^i}^\top] = \delta(t,{t'}) S_{i}$, where $S_{i} \in \mbb{S}^{n_i}_{++}$, and $\delta(t,{t'}) = 1$ if $t={t'}$, and $\delta(t,{t'}) = 0$ otherwise. Additionally, $x_t \in \mbb{R}^n$, $u_t \in \mbb{R}^m$, and $w_t \in \mbb{R}^n$ denote the global state, control input, and process noise of the system, respectively, with $A_t \in \mbb{R}^{n \times n}$ and $B_t \in \mbb{R}^{n \times m}$, where $n \coloneqq \sum_{i=1}^N n_i$ and $m \coloneqq \sum_{i=1}^N m_i$.

Let $\mcl{G} \coloneqq (\mcl{V}, \mcl{E})$ represent the \emph{system graph}, where $\mcl{V} \coloneqq \{1,2,\ldots,N\}$ is a set of vertices (i.e., subsystem indices) and $\mcl{E} \subseteq \mcl{V} \times \mcl{V}$ is a set of edges such that $\forall i,j \in \mcl{V},~ (i,j) \in \mcl{E} \iff \exists t : A_t^{ji} \neq 0$. To capture communication constraints between subsystems, we adopt the following definition \cite{Anderson2019}.

\begin{definition}[$d$-Hop Neighbors]
    For a locality parameter $d \in \mbb{N}$, the sets of \emph{$d$-outgoing neighbors} and \emph{$d$-incoming neighbors} of subsystem $i$ are defined as $\mcl{O}_i (d) \coloneqq \left\{ j \in \mcl{V} : \mathrm{dist} (i,j) \le d \right\}$ and $\mcl{I}_i (d) \coloneqq \left\{j \in \mcl{V} : \mathrm{dist} (j,i) \le d \right\}$, respectively, where $\mathrm{dist}(i,j)$ denotes the distance between (i.e., the number of edges in a shortest directed path connecting) vertices $i,j \in \mcl{V}$.
\end{definition}

The local control policy $\gamma^i$ of subsystem $i$ is said to be \emph{admissible} if each $u_t^i$ depends on $t$ and the past global state histories $\{x_0,x_1,\ldots,x_t\}$. Furthermore, an admissible policy $\gamma^i$ is said to be $d$-\emph{localized} if each $u_t^i$ depends only on $t$ and the past local state histories of its $d$-incoming neighbors:
\begin{align} \label{eqn:localizedcontrollers}
    u_t^i = \gamma^i ( t, \{ x_{0}^j,x_1^j,\ldots,x_t^j \}_{j \in \mcl{I}_i(d)} ).
\end{align}


\begin{problem} \label{prob:mainprob}
    Consider the coupled system \eqref{eqn:linearsystem} with corresponding system graph $\mcl{G}$. Let $\mu_0, \mu_f \in \mbb{R}^{n}$, $\Sigma_0, \Sigma_f \in \mbb{S}_{++}^{n}$, and $d, T \in \mbb{N}$. Let $\{Q_t,R_t\}_{t=0}^{T-1}$ be given such that $Q_t \in \mbb{S}_{+}^{n}$ and $R_t \in \mbb{S}_{++}^{m}$ for all $t$. Assume, without loss of generality, $\mu_f = 0$. Find a set of admissible control policies $\{\gamma^i\}_{i \in \mcl{V}}$ that solves the following localized covariance steering problem:\footnote{Without \eqref{eqn:centralizedlocalityconstraints}, problem \eqref{eqn:centralizedCS} reduces to the standard optimal covariance steering problem, which can be transformed into and solved as a convex semi-definite programming problem; see \cite{bakolas2018finite}.}
    \begin{subequations} \label{eqn:centralizedCS}
        \begin{align}
            & \minimize_{\{\gamma^i\}_{i \in \mcl{V}}} && \sum_{t=0}^{T-1} \mbb{E} \left[ x_t^\top Q_t x_t + u_t^\top R_t u_t \right] \label{eqn:centralizedCScost}
            \\
            & \subjectto  && \eqref{eqn:linearsystem}, \quad \forall t \in \{ 0,1,\ldots,T-1 \},
            \\
            &&& \eqref{eqn:localizedcontrollers}, \quad \forall t \in \{ 0,1,\ldots,T-1 \}, ~ \forall i \in \mcl{V}, \label{eqn:centralizedlocalityconstraints}
            \\
            &&& x_0 \sim \mcl{N}(\mu_0, \Sigma_0),
            \\
            &&&\mbb{E}[x_T] = \mu_f, \label{eqn:finalmeanconstraint}
            \\
            &&&(\Sigma_f - \Cov[x_T]) \in \mbb{S}_{+}^{n}. \label{eqn:finalcovarianceconstraint}
        \end{align}
    \end{subequations}
\end{problem}

Note that the constraint in \eqref{eqn:finalcovarianceconstraint} is a convex relaxation of the non-convex equality constraint, $\Cov [x_T] = \Sigma_f$, which establishes an upper bound on the system's uncertainty in reaching the desired terminal mean $\mu_f$ \cite{bakolas2018finite}.

\begin{remark}
    Existing work on distributed CS overlooks subsystem coupling and instead focuses on decoupled `agents' whose local state evolves according to \cite[Equation (17a)]{Saravanos21}:
    \begin{align} \label{eqn:agent_dynamics}
        x_{t+1}^{i} &= A_{t}^{i} x_{t}^{i} + B_{t}^{i} u_{t}^{i} + w_{t}^{i},
    \end{align}
    with local state matrices $A_{t}^{i} \in \mbb{R}^{n_i \times n_i}$. Unlike the interconnected subsystem dynamics in \eqref{eqn:linearsystem}, each agent's dynamics \eqref{eqn:agent_dynamics} evolves based only on its own state, input, and noise, thereby making it incapable of modeling tightly coupled systems such as power grids.
\end{remark}

\begin{remark} \label{rem:localcovarianceconstraints}
    Additionally, related work imposes terminal covariance constraints of the form \cite[Equation (19)]{Saravanos21}:
    \begin{align} \label{eqn:localconstraints}
        (\Sigma_{f}^{i} - \Cov[x_{T}^{i}]) &\in \mbb{S}_{+}^{n_i}, & \forall i \in \mcl{V},
    \end{align}
    for some $\Sigma_f^i \in \mbb{S}_{++}^{n_i}$, $\forall i \in \mcl{V}$. These `local' constraints ensure that each agent $i$ meets an individual covariance requirement on their own states. In contrast, in Problem \ref{prob:mainprob}, we consider a more general constraint, namely \eqref{eqn:finalcovarianceconstraint}, which acts on the global state of the system. Clearly, this global constraint naturally subsumes the local constraints in \eqref{eqn:localconstraints}.
\end{remark}


\section{Main Results} \label{sec:mainresults}

In this work, we will restrict our attention to admissible control policies associated with causal LTV controllers: $\forall t$, $u_t = K_{t,0} x_{0} + K_{t,1} x_{1} + \ldots + K_{t,t} x_{t}$, where $K_{t,t'} \in \mbb{R}^{m \times n}$, $\forall t' \in \{ 0,1,\ldots,t \}$. One natural way to impose the communication constraints \eqref{eqn:centralizedlocalityconstraints} on the structure of this controller is to enforce the sparsity constraints $[K_{t,t'}]^{ij} = 0$ (for all $t,t'$ and all $i,j$ with $i \notin \mcl{O}_j(d)$), where $[K_{t,t'}]^{ij}$ is a submatrix of $K_{t,t'}$ such that $u_t^i = \sum_{t'=0}^{t} \sum_{j=1}^{N} [K_{t,t'}]^{ij} x_{t'}^j$. However, it is known that for certain (e.g., strongly connected\footnote{That is, the corresponding system graph $\mcl{G}$ is \emph{strongly connected}, meaning there exists a directed path from any vertex to any other vertex in the graph.}) systems, such a constraint can lead to a nonconvex formulation after applying parameterizations \cite[Section~3.5]{Anderson2019}. To circumvent this issue, we will utilize the SLS framework \cite{Anderson2019}, in which the closed-loop response of the system, rather than the controller itself, is directly designed.

\subsection{System Level Synthesis and Locality Constraints} \label{sec:prelim}

Let $\mbf{x} \coloneqq [x_0^\top, x_1^\top, \ldots, x_T^\top]^\top$, $\mbf{u} \coloneqq [u_0^\top, u_1^\top, \ldots, u_{T-1}^\top]^\top$, and $\mbf{w} \coloneqq [x_0^\top, w_0^\top, w_1^\top, \ldots, w_{T-1}^\top]^\top$. Additionally, let $\mbf{A} \coloneqq \begin{bmatrix}
    \blkdiag (A_0, A_1, \dots, A_{T-1}) ~ 0_{Tn \times n}
\end{bmatrix}$, and $\mbf{B} \coloneqq \blkdiag (B_0, B_1, \dots, B_{T-1})$. Using these definitions, the system dynamics \eqref{eqn:linearsystem} can be compactly written as
\begin{align} \label{eqn:closedloopdynamics}
    \mbf{x} = \mbf{Z} \mbf{A} \mbf{x} + \mbf{Z} \mbf{B} \mbf{u} + \mbf{w},
\end{align}
where $\mbf{Z}$ is the block-lower shift matrix, which has identity matrices along its first block subdiagonals and zeros elsewhere. Let $\mbf{K}$ be the block-lower-triangular (BLT) matrix defined by
\begin{align*} \scriptsize
    \mbf{K} \coloneqq
    \begin{bmatrix}
        K_{0,0} & 0 & \cdots & 0 & 0
        \\
        K_{1,0} & K_{1,1} & \cdots & 0 & 0
        \\
        \vdots & \ddots & \ddots & \vdots & \vdots
        \\
        K_{T-1,0} & \cdots & K_{T-1,T-2} & K_{T-1,T-1} & 0
    \end{bmatrix}.
\end{align*}
The closed-loop behavior of the system \eqref{eqn:linearsystem} under the feedback gain $\mbf{K}$ can then be characterized as follows:
\begin{subequations} \label{eqn:systemresponses}
    \begin{align}
        \mbf{x} &= \underbrace{\left(I - \mbf{Z} \left( \mbf{A} + \mbf{B} \mbf{K} \right) \right)^{-1}}_{\eqqcolon \mbf{\Phi}_x} \mbf{w}, \label{eqn:xsystemresponse} \\
        \mbf{u} &= \underbrace{\mbf{K} \left(I - \mbf{Z} \left( \mbf{A} + \mbf{B} \mbf{K} \right) \right)^{-1}}_{\eqqcolon \mbf{\Phi}_u} \mbf{w}, \label{eqn:usystemresponse}
    \end{align}
\end{subequations}
where $\mbf{\Phi}_x$ and $\mbf{\Phi}_u$ are BLT matrices called the state and control \emph{system responses}, respectively. In the sequel, we will write $\Phi_\bullet^{t,t'}$ to denote the $(t,t')$th entry of $\mbf{\Phi}_\bullet$, $\forall \bullet \in \{x,u\}$.

By virtue of the following theorem, instead of optimizing over the feedback gain $\mbf{K}$, the controller synthesis can be performed by optimizing over the system responses $\{\boldsymbol{\Phi}_x$, $\boldsymbol{\Phi}_u\}$.

\begin{theorem}\cite[Theorem 2.1]{Anderson2019} \label{theo:SLS}
    For the LTV system dynamics \eqref{eqn:linearsystem} evolving under the feedback control law $\mbf{u} = \mbf{K} \mbf{x}$, the following statements are true:
    \begin{enumerate}
        \item The affine subspace defined by
        \begin{equation}
            \mbf{Z}_{\mbf{A},\mbf{B}}
            \begin{bmatrix}
                \mbf{\Phi}_x \\ \mbf{\Phi}_u
            \end{bmatrix}
            = I
        \label{eqn:affineSubspace}
        \end{equation}
        parametrizes all possible system responses \eqref{eqn:systemresponses}, where $\mbf{Z}_{\mbf{A},\mbf{B}} \coloneqq \begin{bmatrix}
            I - \mbf{Z} \mbf{A} & -\mbf{Z} \mbf{B}
        \end{bmatrix}$.
        \item For any pair $\{\mbf{\Phi}_x, \mbf{\Phi}_u\}$ of BLT matrices satisfying \eqref{eqn:affineSubspace}, the feedback gain $\mbf{K} = \mbf{\Phi}_u \mbf{\Phi}_{x}^{-1}$ achieves the desired closed-loop response.
    \end{enumerate}
\label{theo:achievableSystemResponses}
\end{theorem}

\begin{remark} \label{rem:identities}
    If constraint \eqref{eqn:affineSubspace} holds, the diagonal entries of $\mbf{\Phi}_x$ (i.e., $\Phi_x^{t,t}$, $\forall t$) are identity matrices, implying $\mbf{\Phi}_x$ is invertible.
\end{remark}

The construction of the feedback gain $\mbf{K}$ requires the inversion of the matrix $\mbf{\Phi}_x$. However, it is known that this can be avoided via the following controller realization (see \cite[Section~2.1]{Anderson2019}):
\begin{align} \label{eqn:controllerrealization}
    u_t = \sum_{t' =1}^{T-1} \Phi_u^{t,t'} \hat{w}_{t-t'},
\end{align}
where $\hat{w}_t = x_{t+1}-\hat{x}_{t+1}$ and $\hat{x}_{t+1} = \sum_{t'=1}^{T-1} \Phi_x^{t+1,t'+1} \hat{w}_{t-t'}$.

The main advantage of the SLS framework is that the communication constraints \eqref{eqn:centralizedlocalityconstraints} can be cast as affine (subspace) constraints on the system responses $\{\mbf{\Phi}_x, \mbf{\Phi}_u\}$, as summarized in the following definitions from \cite{Anderson2019}.\footnote{Note that the $d$-locality constraint \eqref{eqn:localityconstraints} does not necessarily ensure that the feedback gain $\mbf{K} = \mbf{\Phi}_u \mbf{\Phi}_x^{-1}$ itself is $d$-localized. However, the controller realization in \eqref{eqn:controllerrealization} is ensured to be $d$-localized, as any structural constraints imposed on the system responses $\{\mbf{\Phi}_x, \mbf{\Phi}_u\}$ are directly translated; see \cite[Section~2.1]{Anderson2019} or \cite[Section~III-B]{Alonso2023} for further details.}

\begin{definition}[$d$-Localized System Responses]
    Let $[\mbf{\Phi}_{x}]^{ij}$ be the submatrix of $\mbf{\Phi}_x$ that maps the local noise history $\mbf{w}^j \coloneqq [x_0^{j \top},w_0^{j \top},w_1^{j \top},\ldots,w_{T-1}^{j \top}]^\top$ of subsystem $j$ to the local state history $\mbf{x}^i \coloneqq [x_0^{i \top},x_1^{i \top},\ldots,x_{T}^{i \top}]^\top$ of subsystem $i$. Then, $\mbf{\Phi}_x$ is said to be \emph{$d$-localized} if $[\mbf{\Phi}_x]^{ij} = 0$, $\forall j \in \mcl{V}$, $\forall i \notin \mcl{O}_j (d)$. An analogous definition holds for $\mbf{\Phi}_u$.
\end{definition}

\begin{definition}[$d$-Locality Constraints]
    A subspace $\mcl{L}(d)$ is said to be a \emph{$d$-locality constraint} if
    \begin{align} \label{eqn:localityconstraints}
        \begin{bmatrix}
            \mbf{\Phi}_x \\ \mbf{\Phi}_u
        \end{bmatrix} \in \mcl{L}_d
    \end{align}
    implies that $\mbf{\Phi}_x$ is $d$-localized and $\mbf{\Phi}_u$ is $(d+1)$-localized.
\end{definition}

\subsection{Localized SLS Covariance Steering Problem} \label{sec:localizedSLSCSproblem}

We will now employ the SLS framework to transform the localized CS problem \eqref{eqn:centralizedCS} into an SLS problem over system responses. Excluding the constraints \eqref{eqn:centralizedlocalityconstraints} at the moment, problem \eqref{eqn:centralizedCS} can be written as
\begin{subequations} \label{eqn:SLSCS}
\begin{align}
    &\minimize \limits_{\mbf{x}, \mbf{u}} && \mbb{E} \left[ \left\Vert
    \mbf{F}^{\frac{1}{2}}
    \begin{bmatrix}
        \mbf{x} \\
        \mbf{u}
    \end{bmatrix}
    \right\Vert_2^2 \right]
    \\
    &\subjectto && \eqref{eqn:closedloopdynamics}, \\
    &&&\mbf{P}_0 \mbf{x} \sim \mcl{N} (\mu_0, \Sigma_0), \\
    &&&\mbb{E}[\mbf{P}_T \mbf{x}] = \mu_f, \\
    &&&(\Sigma_f - \Cov[\mbf{P}_T \mbf{x}]) \in \mbb{S}_{+}^{n}, \label{eqn:covconst}
\end{align}
\end{subequations}
where $\mbf{Q} \coloneqq \blkdiag (Q_0, Q_1, \dots, Q_{T-1}, 0)$, $\mbf{R} \coloneqq \blkdiag (R_0, R_1, \dots, R_{T-1})$, $\mbf{F} \coloneqq \blkdiag(\mbf{Q},\mbf{R})$, and $\mbf{P}_t$ is a $1 \times (T+1)$ block matrix, where all blocks are $n \times n$ zero matrices except for the $(t+1)$th block which is an identity matrix. After some algebraic manipulations, the mean and covariance of $\mathbf{x}$ can be written in $\mbf{\Phi}_x$ as
\begin{align}
    \mbb{E}[\mbf{x}] = \boldsymbol{\Phi}_x \boldsymbol{\mu}_{\mathbf{w}}, \qquad \Cov[\mathbf{x}] = \boldsymbol{\Phi}_x \boldsymbol{\Sigma}_\mbf{w} \boldsymbol{\Phi}_x^\top, \label{eqn:meancov}
\end{align}
where $\boldsymbol{\mu}_{\mbf{w}} \coloneqq [\mu_0^\top, 0^\top ]^\top$ and
$\boldsymbol{\Sigma}_\mbf{w} \coloneqq \Cov[\mbf{w}] = \blkdiag(\Sigma_0, W_0, W_1, \ldots, W_{T-1})$ with $W_t \coloneqq \Cov[w_t]$. In view of \eqref{eqn:meancov} and Remark \ref{rem:identities}, and with the addition of \eqref{eqn:localityconstraints} in place of \eqref{eqn:centralizedlocalityconstraints}, problem \eqref{eqn:SLSCS} can be formulated as the following convex localized SLS problem:\footnote{Note that each system response is taken from a set of BLT matrices, which is a convex subset of the set of real matrices of the same dimension.}\textsuperscript{,}\footnote{While $d$-locality constraints are always convex, for some $d$, there may not exist system responses satisfying both \eqref{eqn:affineSubspace} and \eqref{eqn:localityconstraints}. Hence, the value of $d$ must be selected with great care (see, for instance, \cite{wang2016acc}). In this work, similar to \cite{Alonso2020DLMPC,Alonso2023}, we assume that, for a given $d$, the localized SLS problem \eqref{eqn:SLScentralizedproblem} is feasible.}
\begin{subequations} \label{eqn:SLScentralizedproblem}
\begin{alignat}{2}
    &\minimize\limits_{\boldsymbol{\Phi}_x, \boldsymbol{\Phi}_u} \quad 
    & &f(\mbf{\Phi}_x,\mbf{\Phi}_u) \coloneqq \left\Vert 
        \mbf{F}^{\frac{1}{2}}
        \begin{bmatrix}
            \mbf{\Phi}_x
            \\
            \mbf{\Phi}_u
        \end{bmatrix}
        \mbf{\Theta}^{\frac{1}{2}}
    \right\Vert_F^2
    \label{eqn:SLScentralizedobjective} \\ 
    &\subjectto \quad 
    & & \eqref{eqn:affineSubspace}, ~ \eqref{eqn:localityconstraints}, \\
    & & &\mbf{P}_T \boldsymbol{\Phi}_x \mbf{P}_0^\top \mu_0 = \mu_f, \label{eqn:SLScentralizedterminal} \\
    & & &\begin{bsmallmatrix}
        \Sigma_f & \mbf{P}_T \boldsymbol{\Phi}_x \boldsymbol{\Sigma}_\mbf{w}^{\frac{1}{2}} \\
        \left(\mbf{P}_T \boldsymbol{\Phi}_x \boldsymbol{\Sigma}_\mbf{w}^{\frac{1}{2}}\right)^\top & I
    \end{bsmallmatrix} \succeq 0, \label{eqn:LMIterminalconstraint}
\end{alignat} 
\end{subequations}
where $\mbf{\Theta} \coloneqq \boldsymbol{\Sigma}_{\mbf{w}} + \boldsymbol{\mu}_{\mbf{w}} \boldsymbol{\mu}_{\mbf{w}}^\top$. The linear-matrix-inequality (LMI) constraint \eqref{eqn:LMIterminalconstraint} is obtained by applying Schur's complement formula to the terminal covariance constraint:
\begin{align*}
    \left(\Sigma_f - \mathbf{P}_T \boldsymbol{\Phi}_x \boldsymbol{\Sigma}_\mbf{w} \boldsymbol{\Phi}_x^\top \mathbf{P}_T^\top \right) \in \mathbb{S}^{n}_{+}.
\end{align*}
For the convexity of constraint \eqref{eqn:LMIterminalconstraint}, see \cite[Proposition~3]{bakolas2018finite}. Note that the localized SLS CS problem \eqref{eqn:SLScentralizedproblem} includes constraint \eqref{eqn:SLScentralizedterminal} that couples the initial and target means and constraint \eqref{eqn:LMIterminalconstraint} that couples the initial and target covariances.

\subsection{Separability of the Localized SLS CS Problem} \label{sec:separability}

Solving the convex but centralized SLS problem \eqref{eqn:SLScentralizedproblem} requires knowledge of all system and problem parameters, including the system matrices $\{\mbf{A},\mbf{B}\}$, cost matrices $\{\mbf{Q},\mbf{R}\}$, and the target mean $\mu_f$ and covariance $\Sigma_f$ of the global system state. However, this assumption can be overly restrictive in practice, as each subsystems $i$ may have access to only local information on the aforementioned variables. To address this scenario, we will investigate the separable structure of problem \eqref{eqn:SLScentralizedproblem} by utilizing the notion of \emph{separability} \cite{wang2018}. This structural insight will be utilized in the subsequent section to distribute computation across subsystems based on their local information and communication constraints.

Let $\mathscr{C} \coloneqq \{\mcl{C}^1,\ldots,\mcl{C}^N\}$ and $\mathscr{R}_\bullet \coloneqq \{\mcl{R}_\bullet^1,\ldots,\mcl{R}_\bullet^N\}$ be partitions of the sets of indices of the columns of $\mbf{\Phi}_x$ (or $\mbf{\Phi}_u$) and the rows of $\mbf{\Phi}_\bullet$, respectively, $\forall \bullet \in \{x,u\}$. Additionally, let $\mathscr{S} \coloneqq \{\mcl{S}^1,\ldots,\mcl{S}^N\}$ be a partition of the set $\{1,\ldots,n\}$.


First, both the system-level parametrization constraint \eqref{eqn:affineSubspace} and the locality constraint \eqref{eqn:localityconstraints} are known to be \emph{column-wise separable} with respect to the column-wise partition $\mathscr{C}$ (see \cite[Section~III-A]{wang2018}): $\eqref{eqn:affineSubspace} \iff
\mbf{Z}_{\mbf{A},\mbf{B}}
\begin{bmatrix}
    \mbf{\Phi}_x(:,\mcl{C}^i) \\ \mbf{\Phi}_u(:,\mcl{C}^i)
\end{bmatrix}
= I(:,\mcl{C}^i)$, $\forall i \in \mcl{V}$; and $\eqref{eqn:localityconstraints} \iff
\begin{bmatrix}
    \mbf{\Phi}_x(:,\mcl{C}^i) \\ \mbf{\Phi}_u(:,\mcl{C}^i)
\end{bmatrix}
\in \mcl{L}_{d}(:,\mcl{C}^i)$, $\forall i \in \mcl{V}$.

Subsequently, the terminal mean constraint \eqref{eqn:SLScentralizedterminal} is \emph{row-wise separable} with respect to the partition $\mathscr{S}$, namely $\eqref{eqn:SLScentralizedterminal} \iff [\mbf{P}_T \mbf{\Phi}_x \mbf{P}_0^\top](\mcl{S}^i,:) \mu_0 = \mu_f(\mcl{S}^i,:), ~ \forall i \in \mcl{V}$.\footnote{Notice that $\mbf{P}_T \boldsymbol{\Phi}_{x} \mbf{P}_0^\top$ is the first block of $\mbf{P}_T \mbf{\Phi}_x$, or equivalently $\Phi_x^{T,0}$.} In contrast, the global terminal covariance constraint \eqref{eqn:LMIterminalconstraint} is generally not separable due to the existence of nondiagonal elements of $\Sigma_f$; as such, exceptions include the special instance with local terminal covariance constraints as described in Remark \ref{rem:localcovarianceconstraints}, in which case $\Sigma_f$ is diagonal.

Lastly, the objective function $f$ in \eqref{eqn:SLScentralizedobjective} is generally nonseparable. However, if $\mbf{F}$ is diagonal, then $f$ is \emph{row-wise separable} with respect to the row-wise partitions $\mathscr{R}_x$ and $\mathscr{R}_u$: $f(\mbf{\Phi}_x,\mbf{\Phi}_u) = \sum_{i \in \mcl{V}} \big( \sum_{r \in \mcl{R}_x^i} \mbf{Q}(r,r) \mbf{\Phi}_{x}(r,:) \mbf{\Theta} \mbf{\Phi}_x(r,:)^\top + \sum_{s \in \mcl{R}_u^i} \mbf{R}(s,s) \mbf{\Phi}_{u}(s,:) \mbf{\Theta} \mbf{\Phi}_u(s,:)^\top \big)$. Further, if $\mbf{\Theta}$ is diagonal, then $f$ is \emph{column-wise separable} with respect to the column-wise partition $\mathscr{C}$: $f(\mbf{\Phi}_x,\mbf{\Phi}_u) = \sum_{i \in \mcl{V}} \sum_{c \in \mcl{C}^i} \big( \mbf{\Theta}(c,c) \mbf{\Phi}_x(:,c)^{\top} \mathbf{Q} \mbf{\Phi}_x(:,c) + \mbf{\Phi}_u(:,c)^{\top} \mathbf{R} \mbf{\Phi}_u(:,c) \big)$. To address general instances in which these conditions do not necessarily hold, the proposition below provides a transformation that will allow one to recast problem \eqref{eqn:SLScentralizedproblem} as an equivalent problem with a separable objective function while maintaining the separable structures of the aforementioned constraints.

\begin{proposition} \label{prop:transformedproblem}
    Let the pair $\{\mbf{\Psi}_x^\star,\mbf{\Psi}_u^\star\}$ of system responses be a solution of the following convex SLS problem:
    \begin{subequations} \label{eqn:transformedSLSproblem}
    \begin{alignat}{2}
        &\minimize \limits_{\mbf{\Psi}_x, \mbf{\Psi}_u} ~~
        & & \tilde{f}(\mbf{\Psi}_x,\mbf{\Psi}_u) \coloneqq \left\Vert 
            \mbf{F}^{\frac{1}{2}}
            \begin{bmatrix}
                \mbf{\Psi}_x
                \\
                \mbf{\Psi}_u
            \end{bmatrix}
            \mbf{\Lambda}^{\frac{1}{2}}
        \right\Vert_F^2 \label{eqn:transformedSLSobjective}
        \\ 
        &\subjectto ~~  
        & & \mbf{Z}_{\mbf{A},\mbf{B}}
            \begin{bmatrix}
                \mbf{\Psi}_x
                \\
                \mbf{\Psi}_u
            \end{bmatrix}
            = \mbf{V}, \label{eqn:transformedSLSparametrization}
        \\
        & & & \begin{bmatrix}
                \mbf{\Psi}_x
                \\
                \mbf{\Psi}_u
            \end{bmatrix} \mbf{V}^\top \in \mcl{L}_d, \label{eqn:transformedlocality}
        \\
        & & & \mbf{P}_T \mbf{\Psi}_x \mbf{V}^\top \mbf{P}_0^\top \mu_0 = \mu_f, \label{eqn:transformedSLSterminal}
        \\
        & & &\begin{bsmallmatrix}
            \Sigma_f & \mbf{P}_T \mbf{\Psi}_x \mbf{V}^\top \bm{\Sigma}_\mbf{w}^{\frac{1}{2}} \\
            \big(\mbf{P}_T \mbf{\Psi}_x \mbf{V}^\top \bm{\Sigma}_\mbf{w}^{\frac{1}{2}}\big)^\top & I
        \end{bsmallmatrix} \succeq 0, \label{eqn:transformedSLSLMI}
    \end{alignat} 
    \end{subequations}
    where $\mbf{\Lambda}$ is a diagonal matrix whose entries are the eigenvalues of $\mbf{\Theta}$ and $\mbf{V}$ is an orthogonal matrix satisfying $\mbf{\Theta} = \mbf{V} \mbf{\Lambda} \mbf{V}^\top$. Then, the pair $\{\mbf{\Phi}_x^\star,\mbf{\Phi}_u^\star\}$ of system responses, where $\mbf{\Phi}_x^\star \coloneqq \mbf{\Psi}_x^\star \mbf{V}^\top$ and $\mbf{\Phi}_u^\star \coloneqq \mbf{\Psi}_u^\star \mbf{V}^\top$, is a solution of problem \eqref{eqn:SLScentralizedproblem}.
\end{proposition}

\begin{proof}
    Since $\mbf{\Theta}$ is a real symmetric matrix, it can always be diagonalized as $\boldsymbol{\Theta} = \mbf{V} \mbf{\Lambda} \mbf{V}^\top$ with some orthogonal matrix $\mbf{V}$. Since $\mbf{V}$ is orthogonal, we have that the transformation $\varphi(\mbf{\Psi}) = \mbf{\Psi} \mbf{V}^\top$ is one-to-one, and that its inverse is given by $\varphi^{-1}(\mbf{\Phi}) = \mbf{\Phi} \mbf{V}$. Under the map $\varphi$, the transformed problem \eqref{eqn:transformedSLSproblem} is equivalent to the original problem \eqref{eqn:SLScentralizedproblem} in the sense that if $\{\mbf{\Phi}_x^\star,\mbf{\Phi}_u^\star\}$ solves problem \eqref{eqn:SLScentralizedproblem}, then $\{\mbf{\Phi}_x^\star \mbf{V},\mbf{\Phi}_u^\star \mbf{V}\}$ solves problem \eqref{eqn:transformedSLSproblem}, whereas if $\{\mbf{\Psi}_x^\star,\mbf{\Psi}_u^\star\}$ solves problem \eqref{eqn:transformedSLSproblem}, then $\{\mbf{\Psi}_x^\star \mbf{V}^\top,\mbf{\Psi}_u^\star \mbf{V}^\top\}$ solves problem \eqref{eqn:SLScentralizedproblem} (see \cite[Section~4.1.3]{boyd_cvxopt}). Note that $\mbf{\Psi}_x^\star$ is invertible as both $\mbf{\Phi}_x^\star$ and $\mbf{V}$ are invertible. It follows from the linearity of $\varphi$ and the orthogonality of $\mbf{V}$ that
    \begin{align*}
        \mbf{K}^\star = \mbf{\Psi}_u^\star {\mbf{\Psi}_{x}^\star}^{-1} = \mbf{\Phi}_u^\star \mbf{V} \mbf{V}^\top {\mbf{\Phi}_{x}^\star}^{-1} = \mbf{\Phi}_u^\star {\mbf{\Phi}_{x}^\star}^{-1}.
    \end{align*}
    This proves the desired result.
\end{proof}

\begin{remark}
    Since $\mbf{\Lambda}$ is diagonal, whether $\mbf{F}$ is diagonal or not, the objective function $\tilde{f}$ in \eqref{eqn:transformedSLSobjective} is column-wise separable with respect to $\mathscr{C}$: $\tilde{f}(\mbf{\Psi}_x,\mbf{\Psi}_u) = \sum_{i \in \mcl{V}} \sum_{c \in \mcl{C}^i} \mbf{\Lambda}(c,c) \big( \mbf{\Psi}_x(:,c)^{\top} \mathbf{Q} \mbf{\Psi}_x(:,c) + \mbf{\Psi}_u(:,c)^{\top} \mathbf{R} \mbf{\Psi}_u(:,c) \big)$. It can be shown that the constraints \eqref{eqn:transformedSLSparametrization} and \eqref{eqn:transformedlocality} remain column-wise separable with respect to $\mathscr{R}_x$ and $\mathscr{R}_u$, while \eqref{eqn:transformedSLSterminal} remains row-wise separable with respect to $\mathscr{S}$. The separability of the LMI constraint \eqref{eqn:transformedSLSLMI} still depends on whether it can be decomposed into local constraints, as described in Remark \ref{rem:localcovarianceconstraints}.
\end{remark}

\subsection{Distributed Solution using Consensus ADMM} \label{sec:admm}

To solve in a distributed manner the localized SLS problem \eqref{eqn:SLScentralizedproblem}, which as shown generally includes nonseparable components, we employ the consensus ADMM framework \cite{Boyd11}. To ensure consensus across subsystems, we let each subsystem $i \in \mcl{V}$ maintain its local copies of system responses $\{ \boldsymbol{\Phi}_x^{(i)}, \boldsymbol{\Phi}_u^{(i)} \}$, while introducing global variables $ \boldsymbol{\Phi}_x^g $ and $ \boldsymbol{\Phi}_u^g $ representing the consensus system responses such that the consensus constraints are: $\boldsymbol{\Phi}_x^{(i)} = \boldsymbol{\Phi}_x^g$ and $\boldsymbol{\Phi}_u^{(i)} = \boldsymbol{\Phi}_u^g$, $\forall i \in \mcl{V}$. The subsystem communication topology is assumed to be determined by the system graph $\mcl{G}$ and a locality parameter $d$. That is, each subsystem can only send information to its $d$-outgoing neighbors and receive information from its $d$-incoming neighbors. To ensure that all subsystems arrive at a consensus for all $d \geq 1$, we make the following assumption.


\begin{assumption} \label{ass:graph}
    The system graph $\mcl{G}$ is strongly connected.
\end{assumption}

\begin{algorithm}[h]
\caption{Distributed and localized covariance steering for coupled systems using consensus ADMM}
\label{alg:cadmm}
\begin{algorithmic}[1]
\State $\mbf{\Phi}_{x,0}^{(i)}, \mbf{\Phi}_{u,0}^{(i)}$, $\mbf{\Omega}_{x, 0}^{(i)}, \mbf{\Omega}_{u, 0}^{(i)}$, and $k \leftarrow 0$ \Comment{Initialization} \label{alg:init}
\While{not converged}
    \For{$i \in \mcl{V}$}{~\textbf{in~parallel}}
        \State Send $\mbf{\Phi}_{x,k}^{(i)}, \mbf{\Phi}_{u,k}^{(i)}$ to neighbors in $\mcl{O}_i (d)$ \label{alg:send}
        \State Receive $\mbf{\Phi}_{x,k}^{(j)}, \mbf{\Phi}_{u,k}^{(j)}$ from neighbors in $\mcl{I}_i (d)$ \label{alg:receive}
        \State $\mbf{\Phi}_{x,k+1}^{(i)}, \mbf{\Phi}_{u,k+1}^{(i)} \gets$ Use \eqref{eqn:primalupdate} \hfill \Comment{Primal update} \label{alg:primal}
        \State $\mbf{\Omega}_{x,k+1}^{(i)}, \mbf{\Omega}_{u,k+1}^{(i)} \gets$ Use \eqref{eqn:dualupdate} \hfill \Comment{Dual update} \label{alg:dual}
    \EndFor
    \State $k \gets k + 1$
\EndWhile
\end{algorithmic}
\end{algorithm}

The proposed consensus ADMM algorithm for the localized SLS problem \eqref{eqn:SLScentralizedproblem} is detailed in Algorithm \ref{alg:cadmm}. First, the local system responses, $\mbf{\Phi}_{x,0}^{(i)}, \mbf{\Phi}_{u,0}^{(i)}$, are initialized to random matrices, while the local dual matrices, $\mbf{\Omega}_{x, 0}^{(i)}, \mbf{\Omega}_{u, 0}^{(i)}$, are set to zero matrices (line \ref{alg:init}). Each subsystem $i$ then sends its local system responses to its $d$-outgoing neighbors (line \ref{alg:send}) and receives local system responses from its $d$-incoming neighbors (line \ref{alg:receive}). Each subsystem then performs primal and dual update procedures (lines \ref{alg:primal} and \ref{alg:dual}), as specified below.

\textbf{Primal Update:} In the case where $\mbf{Q}$ and $\mbf{R}$ are diagonal, each subsystem $i$ updates their local system responses $\mbf{\Phi}_x^{(i)}$ and $\mbf{\Phi}_u^{(i)}$ by solving the following convex SLS subproblem:
\begin{subequations}
\label{eqn:primalupdate}
\begin{align}
    \minimize_{\mbf{\Phi}_x^{(i)}, \mbf{\Phi}_u^{(i)}} ~& f_i(\mbf{\Phi}_x^{(i)},\mbf{\Phi}_u^{(i)}) \nonumber
    \\
    &+\ones^\top (\mbf{\Omega}_{x,k}^{(i)} \circ \mbf{\Phi}_x^{(i)} + \mbf{\Omega}_{u,k}^{(i)} \circ \mbf{\Phi}_u^{(i)})\ones \nonumber \\ & + \rho \sum_{j \in \mcl{I}_i (d)} \left\Vert \mbf{\Phi}_x^{(i)} - \tfrac{1}{2} \left(\mbf{\Phi}_{x,k}^{(i)} + \mbf{\Phi}_{x,k}^{(j)} \right) \right\Vert_F^2 \nonumber
    \\
    &  + \rho \sum_{j \in \mcl{I}_i (d)} \left\Vert \mbf{\Phi}_u^{(i)} - \tfrac{1}{2} \left(\mbf{\Phi}_{u,k}^{(i)} + \mbf{\Phi}_{u,k}^{(j)} \right)     \right\Vert_F^2 \label{eqn:primalobjective}
    \\
    \subjectto ~&
    \mbf{Z}_{\mbf{A},\mbf{B}}
    \begin{bmatrix}
        \mbf{\Phi}_x^{(i)}(:,\mcl{C}^i) \\ \mbf{\Phi}_u^{(i)}(:,\mcl{C}^i)
    \end{bmatrix}
    = I(:,\mcl{C}^i),\\
    &
    \begin{bmatrix}
        \mbf{\Phi}_x^{(i)}(:,\mcl{C}^i) \\ \mbf{\Phi}_u^{(i)}(:,\mcl{C}^i)
    \end{bmatrix}
    \in \mcl{L}_{d}(:,\mcl{C}^i),\\
     &[\mbf{P}_T \mbf{\Phi}_x^{(i)} \mbf{P}_0^\top](\mcl{S}^i,:) \mu_0 = \mu_f(\mcl{S}^i,:), \\
     &\begin{bsmallmatrix}
        \Sigma_f & \mbf{P}_T \boldsymbol{\Phi}_x^{(i)} \boldsymbol{\Sigma}_{\mbf{w}}^{\frac{1}{2}}
        \\
        \left(\mbf{P}_T \boldsymbol{\Phi}_x^{(i)} \boldsymbol{\Sigma}_{\mbf{w}}^{\frac{1}{2}}\right)^\top & I
    \end{bsmallmatrix} \succeq 0, \label{eqn:primal_LMI}
\end{align}
\end{subequations} 
where $\rho > 0$ is a penalty parameter, and
\begin{align*}
    f_i(\mbf{\Phi}_x,\mbf{\Phi}_u) &\coloneqq \sum_{r \in \mcl{R}_x^i}
    \mbf{Q}(r,r) \mbf{\Phi}_x(r,:) \mbf{\Theta} \mbf{\Phi}_x(r,:)^\top
    \\
    &\quad + \sum_{s \in \mcl{R}_u^i} \mbf{R}(s,s) \mbf{\Phi}_u(s,:) \mbf{\Theta} \mbf{\Phi}_u(s,:)^\top.
\end{align*}
For the case of not necessarily diagonal $\mbf{F}$, an analogous convex SLS subproblem can be constructed by decomposing the transformed problem \eqref{eqn:transformedSLSproblem}, which is omitted here due to limited space.

\textbf{Dual Update:} Each subsystem $i$ updates their dual variables according to: $\forall \bullet \in \{x, u\}$,
    \begin{align} \label{eqn:dualupdate}
    \mbf{\Omega}_{\bullet,k+1}^{(i)} &= \mbf{\Omega}_{\bullet,k}^{(i)} + \rho \sum_{j \in \mcl{I}_i (d)} \left( \mbf{\Phi}_{\bullet,k+1}^{(i)} - \mbf{\Phi}_{\bullet,k+1}^{(j)} \right).
\end{align}

Algorithm \ref{alg:cadmm} keeps iterating until all $N$ subsystems reach a consensus, i.e., when the average consensus constraints' residual norms for both system responses reaches a tolerance $\varepsilon > 0$: $\frac{1}{N} \sum_{i \in \mcl{V}} \sum_{j \in \mcl{I}_i (d)} \left\Vert \boldsymbol{\Phi}_\bullet^{(i)} - \boldsymbol{\Phi}_\bullet^{(j)} \right\Vert_F^2 \le \varepsilon$, $\forall \bullet \in \{x, u\}$.

\begin{proposition} \label{prop:convergence}
    Let $\mcl{F}_i(d)$ denote the feasible set of the SLS subproblem \eqref{eqn:primalupdate} for subsystem $i$. If Assumption \ref{ass:graph} holds and, for all $i \in \mcl{V}$, Slater's condition \cite{boyd_cvxopt} holds (i.e., the relative interior of $\mcl{F}_i(d)$ is nonempty), then the sequence of local system responses $(\{\mbf{\Phi}_{x,k}^{(i)},\mbf{\Phi}_{u,k}^{(i)}\})_{k \in \mbb{N}}$ in Algorithm \ref{alg:cadmm} satisfies, as $k \rightarrow \infty$, 1) $\left\Vert \boldsymbol{\Phi}_{\bullet,k}^{(i)} - \boldsymbol{\Phi}_{\bullet,k}^{(j)} \right\Vert_F^2 \rightarrow 0$, $\forall j \in \mcl{I}_i(d)$ and $\forall \bullet \in \{x,u\}$, and 2) $f(\mbf{\Phi}_{x,k}^{(i)},\mbf{\Phi}_{u,k}^{(i)}) \rightarrow p^\star$, where $p^\star$ denotes the optimal value of the localized SLS problem \eqref{eqn:SLScentralizedproblem}.
\end{proposition}

\begin{proof}
    Let $h_i$ be the extended real-valued function defined by:
    \begin{align*}
        h_i(\mbf{\Phi}_x,\mbf{\Phi}_u) \coloneqq
        \begin{cases}
            f_{i}(\mbf{\Phi}_x,\mbf{\Phi}_u)
            , & \textrm{if}~
            \begin{bmatrix}
                \mbf{\Phi}_{x}
                \\
                \mbf{\Phi}_{u}
            \end{bmatrix} \in \mcl{F}_i(d),
            \\
            +\infty, & \textrm{otherwise}.
        \end{cases}
    \end{align*}
    Let $\mathbf{epi} \, h_i \coloneqq \{ (\mbf{\Phi}_x,\mbf{\Phi}_u,\tau) : h_i(\mbf{\Phi}_x,\mbf{\Phi}_u) \leq \tau \}$ be the epigraph of $h_i$. From the definition of $h_i$, it follows that $\mathbf{epi} \, h_i = \{(\mbf{\Phi}_x,\mbf{\Phi}_u,\tau) \in \mcl{F}_i(d) \times \mbb{R} : f_i(\mbf{\Phi}_x,\mbf{\Phi}_u) \leq \tau \}$. Since $\mcl{F}_i(d)$ is a nonempty closed convex set and $f_i$ is a real-valued convex function, $\mathbf{epi} \, h_i$ is a nonempty closed convex set, which is both necessary and sufficient for $h_i$ to be a closed proper convex function \cite[Section~3.2]{Boyd11}. Now, consider the distributed form of the SLS problem \eqref{eqn:SLScentralizedproblem}:
    \begin{subequations} \label{eqn:admm_slscsproblem}
        \begin{align}
            &\minimize_{\{\mbf{\Phi}_x^{(i)},\mbf{\Phi}_u^{(i)}\}_{i \in \mcl{V}}} ~&& \sum_{i \in \mcl{V}} h_i(\mbf{\Phi}_x^{(i)},\mbf{\Phi}_u^{(i)})
            \\
            &\subjectto && \mbf{\Phi}_x^{(i)} = \mbf{\Phi}_x^{(j)}, \quad \forall j \in \mcl{I}_i(d), ~ \forall i \in \mcl{V},
            \\
            &&& \mbf{\Phi}_u^{(i)} = \mbf{\Phi}_u^{(j)}, \quad \forall j \in \mcl{I}_i(d), ~ \forall i \in \mcl{V}.
        \end{align}
    \end{subequations}
    Clearly, problem \eqref{eqn:admm_slscsproblem} is equivalent to its centralized counterpart \eqref{eqn:SLScentralizedproblem} with the same optimal solution and optimal value, as Assumption \ref{ass:graph} ensures that $\mbf{\Phi}_{\bullet}^{(i)} = \mbf{\Phi}_{\bullet}^{(j)}$, $\forall \bullet \in \{x,u\}$, $\forall i,j \in \mcl{V}$. Lastly, Slater's condition implies that strong duality is achieved. In light of all these results, the rest of the proof follows directly from \cite[Appendix~A]{Boyd11}.
\end{proof}

\begin{remark}
    The computational complexity of Algorithm \ref{alg:cadmm} is determined by the primal update procedure (line \ref{alg:primal}), in which each subsystem solves per iteration the SLS subproblem \eqref{eqn:primalupdate} over $O((n^2+nm) T^2)$ variables subject to $O((n+m)T)$ constraints. While this complexity matches that of the centralized SLS problem \eqref{eqn:SLScentralizedproblem} due to the use of the local copies of the global system responses, it is worth highlighting that solving problem \eqref{eqn:primalupdate} only requires each subsystem $i$ to know, besides the global target covariance $\Sigma_f$, its assigned portion of the system matrices $\{\mbf{A},\mbf{B}\}$, the cost matrices $\{\mbf{Q},\mbf{R}\}$, and the target mean $\mu_f$. Additionally, in the special instance of local terminal covariance constraints as described in Remark \ref{rem:localcovarianceconstraints}, this complexity can further be reduced.
\end{remark}


\section{Numerical Simulations} \label{sec:numericalsimulations}

\begin{figure}[t]
    \centering
    \includegraphics[width=0.22\textwidth]{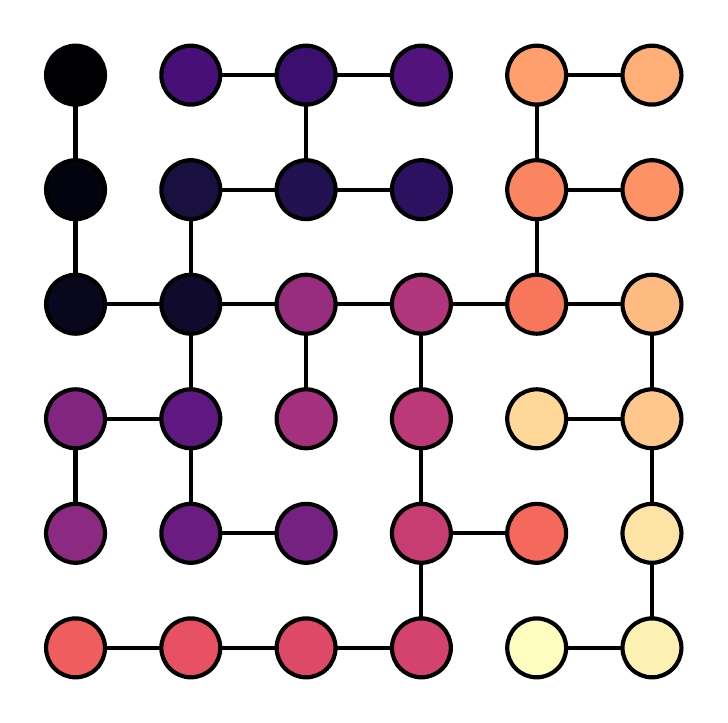}
    \caption{Network topology of the $6 \times 6$ power grid.}
    \label{fig:network_topology}
\end{figure}

In this section, we verify the performance of Algorithm \ref{alg:cadmm} on a power system modeled as a randomized spanning tree within a $6 \times 6$ grid (i.e., $N=36$), as shown in Figure \ref{fig:network_topology}. Each subsystem is governed by the discretized swing equations: $m^i \ddot{\theta}^i + d^i \dot{\theta}^i = - \sum_{j \in \mcl{N}_i} k^{ij} (\theta^i - \theta^j) + w^i + u^i$, where $\theta^i$ and $\dot{\theta}^i$ denote the phase angle and frequency deviations; $m^i$ and $d^i$ are the inertia and damping; $w^i$ and $u^i$ represent external disturbances and control inputs; $k^{ij}$ is the coupling term between subsystems $i$ and $j$; and $\mcl{N}_i$ is the set of neighboring vertices. Defining the local state vector of each subsystem as $x^i \coloneqq \begin{bmatrix}\theta^i & \dot{\theta}^i \end{bmatrix}^\top$, the discretized swing dynamics can be expressed as
\begin{equation}
    x_{t+1}^{i} = A^{ii} x_{t}^{i} + \sum_{j \in \mcl{N}_{i} \backslash \{i\}} A^{ij} x_{t}^{j} + B^{i} u_{t}^{i} + w_{t}^{i},
\end{equation}
where $A^{ii} = \begin{bmatrix}
    1 & \Delta t \\ - \frac{k^i}{m^i} \Delta t & 1 - \frac{d^i}{m^i} \Delta t
\end{bmatrix}$, $
A^{ij} = \begin{bmatrix}
    0 & 0 \\ \frac{k^{ij}}{m^i} \Delta t & 0
\end{bmatrix}$,
$B^{i} = \begin{bmatrix} 1 
& 0 \end{bmatrix}^\top$, $k^i = \sum_{j\in \mcl{N}_i} k^{ij}$, and $\Delta t = 0.2$. The values of $k^{ij}$, $d^i$, and $m^{i}$ are sampled uniformly at random from $[0.5, 1]$, $[0.2, 0.8]$, and $[0.5, 1]$, respectively. Note that the system graph $\mcl{G}$ coincides with the network topology displayed in Figure \ref{fig:network_topology}.

The localized CS problem \eqref{eqn:centralizedCS} is initialized with $T = 10$, $Q = \diag(100,500,\ldots,100,500)$, $R = 0.01 I$, $W_t = 0.2I$, $\mu_f = 0$, and $\Sigma_f = MM^\top$ with the diagonal (resp., nondiagonal) elements of $M$ sampled from $\mcl{N}(0.5, 0.1)$ (resp., $\mcl{N}(0, 0.1)$). The initial state is sampled as $x_0 \sim \mcl{N}(\mu_0, \Sigma_0)$, where $\mu_0$ follows the standard normal distribution scaled by $30$, and $\Sigma_0$ is a diagonal matrix with entries sampled from a uniform distribution over $(0, 60]$. The locality parameter is chosen as $d = 1$, i.e., the communication graph of subsystems is also equivalent to the system graph $\mcl{G}$. We execute Algorithm \ref{alg:cadmm} using \textsc{Clarabel} \cite{Clarabel_2024}, an interior point conic solver available in \textsc{CVXPY} \cite{diamond2016cvxpy}. With $\rho=0.01$ and $\varepsilon = 10^{-4}$, the algorithm converges within approximately $1,000$ iterations.

\begin{figure}[t]
    \centering
    \begin{subfigure}[b]{0.24\textwidth}
        \centering
        \includegraphics[width=\textwidth]{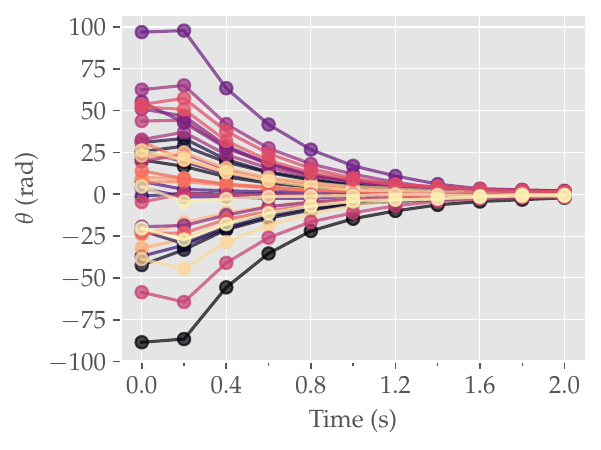}
        \caption{Phase angle deviation vs. time}
    \label{fig:theta}
    \end{subfigure}
    \begin{subfigure}[b]{0.24\textwidth}
        \centering
        \includegraphics[width=\textwidth]{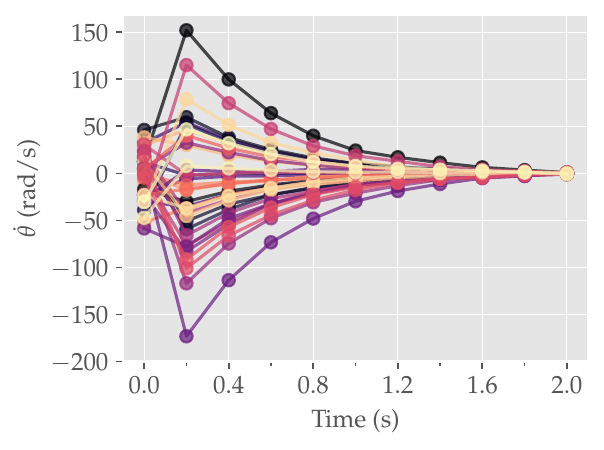}
        \caption{Frequency deviation vs. time}
        \label{fig:theta_dot}
    \end{subfigure}
    \begin{subfigure}[b]{0.24\textwidth}
        \centering
        \includegraphics[width=\textwidth]{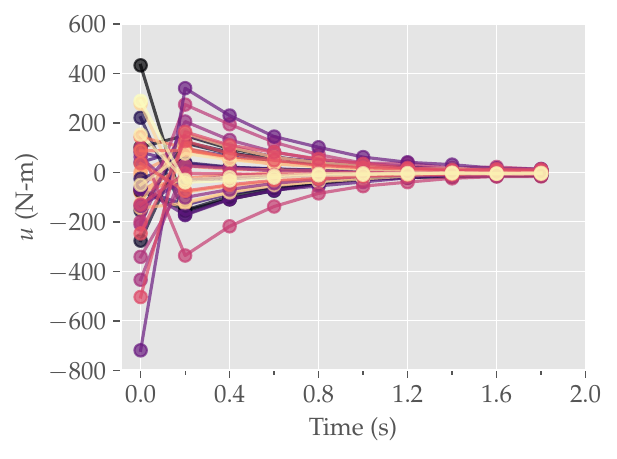}
        \caption{Torque input vs. time}
        \label{fig:torque}
    \end{subfigure}
    \begin{subfigure}[b]{0.24\textwidth}
        \centering
        \includegraphics[width=\textwidth]{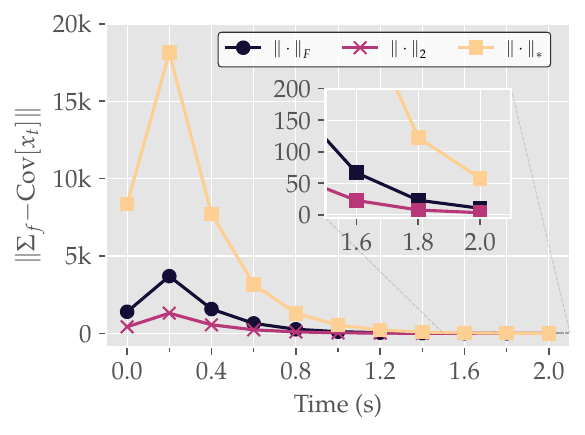}
        \caption{Covariance deviation vs. time}
        \label{fig:norm_difference}
    \end{subfigure}
    
    \caption{Results for the simulation example with $d = 1$.}
    \label{fig:results}
\end{figure}

Figure \eqref{fig:results} displays results for the numerical simulation. Figures \ref{fig:theta} and \ref{fig:theta_dot} display each subsystem's phase angle and frequency deviations over the time horizon, respectively. It can be seen that both state deviations converge to 0 at the final time. Figure \ref{fig:torque} shows the torque input at each time step. Figure \ref{fig:norm_difference} shows the difference between the final covariance and the state covariance at each time step, evaluated using three different norms\footnote{Here, $\|\cdot\|_2$ and $\|\cdot\|_*$ denote the spectral norm and the nuclear (or trace) norm of a matrix, respectively.}. In all cases, the state covariance at the final stage either converges to or closely approximates the target covariance.

\section{Conclusion} \label{sec:conclusion}

This work addressed the problem of localized covariance steering in coupled stochastic linear systems. We employed the system-level synthesis framework to transform this problem into a convex optimization problem over system responses. A consensus-based algorithm was then proposed to distribute computation among subsystems. The effectiveness of this approach was demonstrated on a power system. In future work, we plan to address how clustering agents can accelerate the convergence of the proposed algorithm, as in \cite{lee2025}, and to investigate potential data-driven extensions of distributed and localized covariance steering.

\section*{Acknowledgment}
The first author would like to thank Dr. Takashi Tanaka for the insightful discussions that have led to this work. The authors would also like to thank the editors and anonymous reviewers for their constructive reviews and recommendations for potential extensions.

\bibliographystyle{unsrt}
\bibliography{root}

\addtolength{\textheight}{-12cm}

\end{document}